\documentclass[sn-mathphys]{sn-jnl}

\jyear{2021}%

\theoremstyle{thmstyleone}%

\theoremstyle{thmstyletwo}%

\theoremstyle{thmstylethree}%

\graphicspath{{./figures}}
\usepackage{array,amsmath,amsfonts,amssymb,amsbsy,latexsym}
\usepackage{soul}

\begin{document}

\title[A view of computational models for image segmentation]{A view of computational models for image segmentation}


\author[1]{\fnm{Laura} \sur{Antonelli}}\email{laura.antonelli@cnr.it}
\equalcont{The authors contributed equally to this work.}

\author[2]{\fnm{Valentina} \sur{De Simone}}\email{valentina.desimone@unicampania.it}
\equalcont{The authors contributed equally to this work.}

\author*[3]{\fnm{Daniela} \sur{di Serafino}}\email{daniela.diserafino@unina.it}
\equalcont{The authors contributed equally to this work.}

\affil[1]{\orgdiv{Institute for High Performance Computing and Networking (ICAR)}, \orgname{CNR}, \orgaddress{\street{Via Pietro Castellino, 111}, \city{Napoli}, \postcode{80131}, \state{Italy}}}

\affil[2]{\orgdiv{Department of Mathematics and Physics}, \orgname{University of Campania ``Luigi Vanvitelli''}, \orgaddress{\street{viale Abramo Lincoln, 5}, \city{Caserta}, \postcode{81100}, \state{Italy}}}

\affil*[3]{\orgdiv{Department of Mathematics and Applications ``R.~Caccioppoli''}, \orgname{University of Naples Federico II}, \orgaddress{\street{Via Cintia, Monte S. Angelo}, \city{Napoli}, \postcode{80126}, \state{Italy}}}


\abstract{Image segmentation is a central topic in image processing and computer vision and a key issue in many applications, e.g., in medical imaging, microscopy, document analysis and remote sensing. According to the human perception, image segmentation is the process of dividing an image into non-overlapping regions. These regions, which may correspond, e.g., to different objects, are fundamental for the correct interpretation and classification of the scene represented by the image. The division into regions is not unique, but it depends on the application, i.e., it must be driven by the final goal of the segmentation and hence by the most significant features with respect to that goal. Thus, image segmentation can be regarded as a strongly ill-posed problem. A classical approach to deal with ill posedness consists in incorporating in the model a-priori information about the solution, e.g., in the form of penalty terms.
In this work we provide a brief overview of basic computational models for image segmentation,
focusing on edge-based and region-based variational models, as well as on statistical and machine-learning approaches. We also sketch numerical methods that are applied in computing solutions to these models. In our opinion,
our view can help the readers identify suitable classes of methods for solving their specific problems.}




\maketitle

\section{Introduction}\label{sec:intro}

Image segmentation is a fundamental task of image processing, image analysis, image understanding, and pattern recognition. It has a very long history, whose origin may be dated back to about 50 years ago. A seminal paper is~\cite{bib:BriceFennema1970}, where the authors pointed out that an important component of the Stanford Research Institute automation project was a set of programs providing the automaton with a means of interpreting visual data.

\begin{figure}[ht!]
\medskip
\begin{center}
\newcolumntype{C}{>{\centering\arraybackslash} m{.22\textwidth} }
\begin{tabular}{CCC}
\includegraphics[width=.22\textwidth]{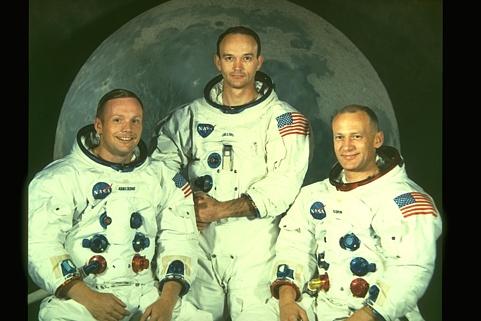} &
\includegraphics[width=.22\textwidth]{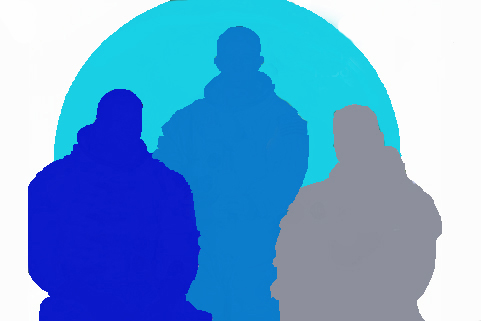} &
\includegraphics[width=.22\textwidth]{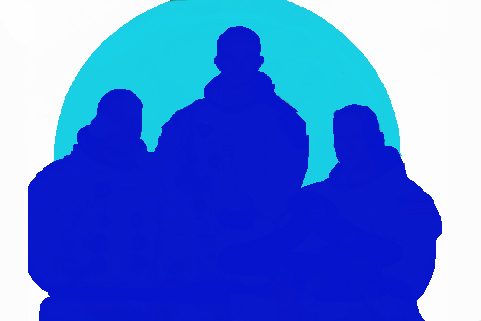}\\[1mm]
\multicolumn{1}{c}{\small Image \texttt{\#323016}}  & \multicolumn{1}{c}{\small \texttt{instance segmentation}}  &
\multicolumn{1}{c}{\small \texttt{semantic segmentation}} \\
& \multicolumn{1}{c}{\small \texttt{(a)}}  &
\multicolumn{1}{c}{\small \texttt{(b)}} \\
\end{tabular}
\caption{Illustration of instance and semantic segmentation of the Berkeley database image \texttt{\#323016}. The results displayed in \texttt{(a)} and \texttt{(b)} were produced by using Adobe Photoshop.\label{fig:segclassification}}
\end{center}
\end{figure}

While it is possible to accurately represent the information in a real scene by an image, this representation alone does not enable us to highlight specific properties of the scene. Conversely, a description in terms of ``natural'' elements of the image, such as regions and boundaries of the visualized objects, represented in a uniform manner, provides easy access to useful global information, thus allowing recognition and extraction of specific image features. Thus, to generate a description of specific elements of the image, it is customary to \emph{segment} the image into more parts (or \emph{segments}). Figure~\ref{fig:segclassification} shows an example of two main types of segmentation, i.e., \emph{instance segmentation}, which identifies the object instance of each pixel for every known object within an image, and \emph{semantic segmentation}, which identifies the object category of each pixel for every known object within an image.

Image segmentation is used in many application fields, such as medical imaging~\cite{bib:Khalid2020Survey}, microscopy imaging~\cite{bib:Bui2020Microscopy}, remote sensing~\cite{bib:Hossain2019RS},
and document image analysis~\cite{bib:Eskenazi2017SegDOC}.
The choice between semantic and instance segmentation is generally dependent on the goal of the classification or object detection step that follows the segmentation phase. For example, in the segmentation of terrain in satellite imagery, we may use the semantic segmentation to distinguish different land areas, like vegetation, ground, water and building, or we may use the instance segmentation to distinguish different common weeds in agricultural fields (i.e., separate instances of objects belonging to the same class).

\begin{figure}[!ht]
\medskip
\begin{center}
\newcolumntype{C}{>{\centering\arraybackslash} m{.22\textwidth} }
\begin{tabular}{CCC}
\includegraphics[width=.22\textwidth]{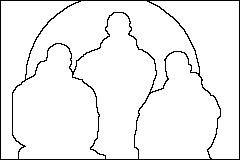} &
\includegraphics[width=.22\textwidth]{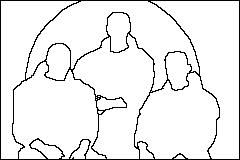} &
\includegraphics[width=.22\textwidth]{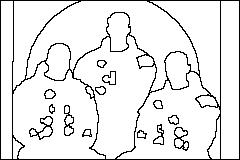} \\[1mm]
\multicolumn{1}{c}{\small User \texttt{\#1103} - 6 segments}  &
\multicolumn{1}{c}{\small  User \texttt{\#1105} - 18 segments}  &
\multicolumn{1}{c}{\small User \texttt{\#1106} - 36 segments}
\end{tabular}
\caption{Segmentations of the Berkeley database image \texttt{\#323016} by different users
        (\small{see \url{https://www2.eecs.berkeley.edu/Research/Projects/CS/vision/grouping/resources.html}}.)
\label{fig:3segastronauti}}
\end{center}
\end{figure}

Since different applications 
may require different partitions to extract significant features, there is no single standard method for image segmentation. Thus, the segmentation problem has not a unique result, as shown in Figure~\ref{fig:3segastronauti}, where different segmentations of the same image are shown, resulting from different segmentation criteria. On the other hand, different methods are not equally effective in segmenting a specific type of image (real scenes, synthetic images, medical images, etc.), and the criteria to define a successful segmentation depend on the desired goal of the segmentation itself.  Therefore, segmentation
remains a challenging problem in image processing and computer vision, in spite of several decades of research.
%
%

We present image segmentation as a highly ill-posed problem, and discuss basic models that take into account a-priori information about the solution, attempting to put these models into a coherent mathematical framework. We look at the inclusion of a-priori information as a sort of regularization approach and show that it is ubiquitous in image segmentation models, from older ``classical'' ones to machine learning approaches, revealing links and similarities between them. Note that we focus on basic models in order to keep our discussion easy and get rid of technical details. We also sketch some numerical methods used in the application of the various models. Although this is only a simplified and partial view of image segmentation, we believe that it may give a contribution towards a better understanding of this huge field.

The rest of this paper is organized as follows. 
In Section~\ref{sec:formulation} we present a mathematical formulation of image segmentation, and in Section~\ref{sec:basic_segm} we discuss basic segmentation models, focusing on edge-based, region-based, statistical and machine learning ones. In Section~\ref{sec:methods} we give a quick overview of numerical techniques that may be used to solve the aforementioned models. Finally, we give some conclusions in Section~\ref{sec:conclusions}.


\section{Mathematical formulation of image segmentation\label{sec:formulation}}

Let $\mathcal{I}$ be the set of the images defined in a domain $\Omega \subset \mathbb{R}^d$ ($d \geq 2$), $I_0 \in \mathcal{I}$ the observed image, and $\mathcal{P}_1, \ldots, \mathcal{P}_n$ logical predicates used to check $n$ statements, expressed using features of the image, e.g., edges, smoothness, texture or color, so that $\mathcal{P}_k(A)=true$ if all the points of $A\subseteq \Omega$ satisfy the $k$-th statement. Just to give an example, in order to compute a two-region segmentation of a normalized gray-level image $I_0$, we can define $n=2$ simple statements as follows, which involve the gray level of the intensity light to separate the background from the foreground:
\[
\mathcal{P}_1(A)= \left\{  I^*(x) < \alpha, \forall x \in A \right\}
, \quad
\mathcal{P}_2(A)= \left\{  I^*(x) \geq \alpha, \forall x \in A \right\}
\]
\noindent
where $I^*$ is a suitable approximation of $I_0$ and $\alpha \in (0,1)$ is a suitable value.

Generalizing the definition in~\cite{bib:Pal1993Review}, the instance segmentation $S$ of $I_0$ according to the predicates $\mathcal{P}_k$, $k=1, \ldots, n$, consists of finding a  decomposition of $\Omega$ into $m$ components $\Omega_i$, with $i = 1, \ldots, m$ and $m \geq n$, such that
\begin{enumerate}
    \item  $\Omega_i \not = \emptyset, \;  \forall \, i \in \{ 1, \ldots, m \}$;
    \item  $\overset{\circ}{\Omega_i} \bigcap \overset{\circ}{\Omega_j} = \emptyset , \; \forall \, i,j \in \{ 1, \ldots, m \}$ with $i \not = j$, where $\overset{\circ}{\Omega_k}$ denotes the interior of $\Omega_k$;
    \item  $\bigcup\limits_{i=1}^{m} \Omega_i = \Omega$;
    \item $\forall i \in \{ 1, \ldots, m \} \; \exists ! \; k \in \{ 1,...,n \} \;$  such that 
    \begin{itemize}
        \item[i.] $\mathcal{P}_k(\Omega_i) =$ \emph{true};
        \item[ii.] $\mathcal{P}_k (\Omega_i \bigcup \Omega_j) =$ \emph{false}, $\forall \, j \in \{ 1, \ldots, m \}$ with $j \not = i$.
    \end{itemize}
\end{enumerate}
By adding to item 4
\begin{itemize} \setlength\itemindent{13pt}
    \item[iii.] $\mathcal{P}_k(\Omega_j) =$ \emph{false},  $\forall \, j \in \{ 1, \ldots, m \}$  with $j \not = i$.
\end{itemize}
we also obtain the semantic segmentation.

We can define the segmentation $S$ of $I_0$ as follows too. Let $\Sigma$ be the set of possible segmentations of the images in $\mathcal{I}$ according to some criteria defined by the predicates $\mathcal{P}_k$.
Then $S$ can be expressed as
\begin{equation*}
S = (u^*,I^*),
\end{equation*}
\noindent
where $u^*$ is a curve that matches the boundaries of the decomposition of $\Omega$, i.e., $u^*= \cup_i \partial \Omega_i$\footnote{The parametric representation of the curve $u$ is defined by a continuous map $\gamma: X \longrightarrow \mathbb{R}^d$, where $X \subset \mathbb{R}$ is an interval and $u=\gamma(X)$. With a little abuse of notation we identify the curve $u$ with the function $\gamma$.}, and $I^*$ is a piecewise-smooth function defined on $\Omega$ that approximates~$I_0$. In particular, we may assume that the restriction of $I^*$ to any set $\overset{\circ}{\Omega_i}$ is piecewise differentiable.
The segmentation $S$ may be also identified directly by using a labeling operator~$\Phi$, i.e.,
\begin{equation}\label{eq:labeling}
S=\Phi(I^*),
\end{equation}
\noindent where
\[ \Phi(I(x)) = l_i \mbox{ if } x \in \Omega_i, \]
$I(x)$ is the value of $I$ associated with $x$, and $l_i \in \, \mathcal{N} = \{ l_1, l_2, \ldots, l_m \}$ is a label.


\section{Basic segmentation models\label{sec:basic_segm}}

We look at image segmentation as an ill-posed problem, whose solution is highly undetermined.
Classical approaches for computing a solution of an ill-posed problem require additional information that enforces uniqueness and stability. To this end, suitably defined penalty terms can be applied. Then, the solution is obtained by minimizing an energy functional~$E$ containing a fidelity term $\mathcal{F}$ that measures the consistency of the candidate segmentation with the observed image, and a penalty term $\mathcal{P}$ that promotes solutions with suitable properties:
\begin{equation}\label{eq:generalE}
    (I^*,u^*) := \arg \, \underset{(I,u)}{\min} \, E(I, u; I_0) = \arg \, \underset{(I,u)}{\min} \left( \mathcal{F}(I,u;I_0) +  \lambda \mathcal{P} (I,u) \right).
\end{equation}
\noindent
Here $\lambda > 0$ is a parameter that generally needs careful tuning to suitably balance $\mathcal{F}$ and $\mathcal{P}$ (see, e.g., \cite{bib:Antonelli2020Adaptive} and the references therein).

The minimization problem~\eqref{eq:generalE} can be solved by writing the Euler-Lagrange equations, which can be derived by integrating by parts the energy functional and using the Gauss theorem along with the fundamental lemma of the calculus of variations. Then a numerical solution can be computed by applying a gradient descent approach, where the descent direction is parameterized through an artificial time, and by a finite-difference discretization. A widely used and effective alternative consists in discretizing problem~\eqref{eq:generalE} and then solving it by a numerical optimization method. We will come back to these two approaches in Section~\ref{sec:methods}.

Recently, machine learning techniques have been successfully applied to segmentation problems. The key idea is to tune a generic model to a specific solution through learning against sample data (training data). The learning phase extracts prior information to be embedded into the penalty term from a large dataset containing pairs of type \emph{(image, ground-truth label)}~\cite{bib:Lucasetal18}. Machine learning approaches using unlabeled image data as training datasets are also available.
Although these techniques successfully solve image segmentation problems, sometimes outperforming state-of-the-art variational models, they have been often designed on-demand for specific tasks used as ``black-box'' models and require a large amount of data to produce results.

In the next subsections we provide some examples of image segmentation models. Note that we focus on basic models, with the aim of providing a general idea of these approaches while avoiding technical details that are outside the scope of this work. It is also worth observing that these models are the basis of modern ones, developed either to improve the effectiveness of the original models in some applications~\cite{bib:CTETRIS22} 
or to complement and refine Machine Learning techniques for segmentation 
~\cite{bib:YousefiriziPET-CT}.

%


\subsection{Edge-based models\label{subsec:edge-based}}

Edge-based models aim at finding $u^*= \cup_i \partial \Omega_i$ by solving the minimization problem~\eqref{eq:generalE} with respect to the curve $u$ (note that $I$ and $I^*$ are not explicitly considered in this case). These models include
the so-called \emph{Active Contours}~\cite{bib:Kasseetall} or \emph{Snakes}. Here the fidelity and regularization terms act as an internal force and an external force, respectively, which move the curve within the image to find the boundaries of the sets $\Omega_i$. More precisely, the energy functional takes the form
\begin{equation}\label{eq:ACmodel}
  E_{AC}(u) = \underbrace{\int_0^1 g(\vert \nabla I_0(u(s))\vert)^2 ds}_{\mathcal{F}} + \lambda \underbrace{\int_0^1 \vert u'(s) \vert^2 ds}_\mathcal{P} ,
\end{equation}
\noindent
where $I_0$ is the observed image, $g$ is an edge-detector function and the curve $u$ is parametrized by $s \in [0,1]$. The first term attracts the curve toward the boundaries, whereas the second one controls its smoothness, and as a result the curve $u$ changes its shape like a snake.

The evolving curve is driven by surface properties, such as curvature and normal direction, and by image features, such as gray levels, hue or saturation in color images, and intensity gradient in 2D images or change in slope in 3D ones. For example, the mean curvature can be used and in this case the edge-detector function is also responsible for stopping the curve on the edges. The function $g$ may be defined as
\[ g(\vert \nabla I_0 \vert) = \frac{1}{1+\vert\nabla (G_\sigma * I_0)\vert^2}, \]
\noindent
where $g$ is a positive and decreasing function, $G_\sigma$ is the Gaussian kernel with standard deviation $\sigma$, and $*$ denotes the convolution operator.

In a Lagrangian approach, an initial curve is evolved by
\begin{equation}\label{pde}
     \frac{\partial u}{\partial t}+ \mathcal{L}(u)=0,
\end{equation}
\noindent
where $\mathcal{L}$ is a differential operator. The simplest evolution is given by $\mathcal{L}(u)= F N$, where $N$ is the normal to the curve and $F$ is a constant that determines the speed of evolution. More generally, the evolution is driven by an external force. For example, in the mean-curvature evolution, $\mathcal{L}(u) = \kappa N$, where $\kappa$ is the Euclidean curvature of~$u$~\cite{bib:AlvarezMorel}.

When $u$ has an explicit representation, it is not easy to deal with topological changes like merge and split, and a re-parametrization of the curve may be required. 
Therefore, the evolution of the curve $u$ is commonly described by level-set methods~\cite{bib:OsherSethian}, thanks to their ability to follow topology changes, cusps and corners. In a level set approach, the curve $u$ is implicitly represented by the zero-level set of a function $\phi(t,x)$, i.e.,
$u = \{ x\in \Omega : \phi(t,x) = 0 \}$. The level set formulations of the simplest evolution and the mean-curvature one read, respectively:
\[
     \frac{\partial \phi}{\partial t} = F \vert \nabla \phi \vert, \;F \in \mathbb{R} \;\;\; \mbox{ and } \;\;\;
     \frac{\partial \phi}{\partial t} = \mbox{div} \left( \frac{\nabla \phi}{\vert \nabla \phi \vert} \right) \vert \nabla \phi \vert.
\]

\subsection{Region-based models}

Region-based models provide directly the segmentation by means of the image partition $\{\Omega_i, \; i=1, \ldots, m \}$. Region-growing models are among the simplest models falling in this class, 
and
in order to obtain accurate segmentations they have been merged with variational approaches where the evolution changes according to the minimization of an energy functional including region-based terms~\cite{bib:RegGrowing}.

A very popular region-growing model was proposed by \emph{Mumford and Shah}~\cite{bib:MumfordShah1989}. In this case, the functional $E$ in~\eqref{eq:generalE} takes the form
\begin{equation} \label{eq:MSmodel}
    E_{MS}(I,u) = \underbrace{\int_\Omega (I-I_0)^2 d x}_\mathcal{F} + \underbrace{\lambda \int_{\Omega - u} \vert \nabla I \vert^2 d x +  \mu \, len(u)}_\mathcal{P},
\end{equation}
\noindent
where $len(u)$ denotes the length of $u$, and $\lambda$ and $\mu$ are positive parameters. The term $\mathcal{F}$ attempts to achieve the minimum distance between $I_0$ and its piecewise-smooth approximation $I$, and $\mathcal{P}$ attempts to reduce the variation of $I$ within each set $\Omega_i$ while keeping the curve $u$ as short as possible. Minimizing~\eqref{eq:MSmodel} in a suitable space provides an optimal pair $(I^*, u^*)$ representing a simplified description of $I_0$ by means of a function with bounded variation and a set of edges~\cite{bib:MumfordShah1989}. Finally, in~\cite{bib:GemanG} the Mumford and Shah model is formulated as a deterministic refinement of a probabilistic model for image restoration.

A simplified version of the Mumford-Shah model is its restriction to piecewise-constant functions. The \emph{Chan-Vese} model~\cite{bib:ChanVese2001} is a particular case of the simplified version, aimed at obtaining a two-phase segmentation where the piecewise-constant function assumes only two values. Its functional $E$ takes the following form:
\begin{equation}\label{eq:CVmodel}
  \begin{array}{rl}
    \displaystyle E_{CV} (I, c_{in},c_{out}) = & \displaystyle
              \!\! \underbrace{\left( \int_\Omega H(I) \left( c_{in} - I_0 \right)^2dx + \int_\Omega (1-H(I)) \left( c_{out} - I_0 \right)^2dx \right)}_\mathcal{F} \\
                                                                 + & \displaystyle \!\! \lambda \, \underbrace{\int_\Omega \vert \nabla H(I)\vert \, dx}_{\mathcal{P}} ,
  \end{array}
\end{equation}

\noindent
where $H$ is the Heaviside function and $c_{in}$ and $c_{out}$ are the average values of the intensity in the foreground and background of the image, respectively. The solution $I^*$ is the best approximation to $I_0$ among all the functions that take only two values.

Minimizing~\eqref{eq:CVmodel} is a nonconvex problem, thus solution methods may get stuck into local minima and result in unsatisfactory segmentations. Aiming to overcome this drawback, some strategies have been proposed, including the convexification of the functional by taking advantage of its geometric properties. An example is given by the two-phase partitioning model introduced by~\emph{Chan, Esedo\=glu and Nikolova}~\cite{bib:CEN2006}:

\begin{equation} \label{eq:CENmodel}
 E_{CEN}(I, c_{in}, c_{out}) =  \underbrace{\int_{\Omega} \left( (c_{in} - I_0 )^2 I + (c_{out} - I_0 )^2  (1-I) \right) dx}_{\mathcal{F}}
                        + \lambda \, \underbrace{\int_{\Omega} \vert \nabla I \vert \, dx}_{\mathcal{P}}
\end{equation}
\noindent
with  $0 \le I \le 1$ and $ c_{in},\, c_{out} > 0$.


\subsection{Statistical models}

Statistical models usually provide a conditional probability, $P(S \vert I_0)$, of a segmentation $S \in \Sigma$ given the observed image $I_0$, and then select
the segmentation with the highest probability. In the \emph{Maximum a Posteriori (MAP)} approach the segmentation is given by
\[  S^* = \arg \, \underset{S \in \Sigma}{\max} \, P(S \vert I_0).  \]
\noindent
According to the Bayes rule,
\[ P(S \vert I_0) = \frac{P(I_0 \vert S) P(S)}{P(I_0)}, \]
where $P(S)$ is the prior probability measuring how well $S$ satisfies certain properties of the given image, and $P(I_0 \vert S)$ is the conditional probability measuring the likelihood of $I_0$ given $S$ (see, e.g., \cite{bib:CalvettiSomersalo2018}).
Since the probability $P(I_0)$ is constant, the segmentation can be obtained by maximizing $P(I_0 \vert S) P(S)$.

Markov Random Field (MRF) models offer a framework to define prior and likelihood by capturing properties of the image such as texture, color, etc.~\cite{bib:Kato}. The segmentation is formulated within an image labeling framework, i.e., \mbox{$S = \Phi(I(x))$}, where the problem is reduced to find the labeling which maximizes the posterior probability. Label dependencies are modeled by an MRF. Then, using  the Hammersley-Clifford theorem, we get the Gibbs distribution
\[ P(S) = \frac{1}{Z}exp(-U(S)), \]
\noindent
where the energy function $U$ takes the form
\[ U(S) = \sum_{c \in C} V_c(S_c), \]
$C$ is the set of cliques of $S$, $V_c(S_c)$ is the potential of the clique $c \in C$ having the label configuration $S_c$, and $Z$ is a normalizing constant. The conditional probability $P(I_0 \vert S)$ can be modeled by a Gaussian distribution. Then the original MAP estimation is equivalent to the following energy minimization problem:
\[  S^* = \arg \, \underset{S}{\max} \, P(S \vert I_0) =  \arg \, \underset{S}{\min}  \, U(S). \]


\subsection{Machine learning models}

Machine learning approaches, and in particular deep learning ones, are more and more used in solving image segmentation problems, also outperforming the previous approaches. Roughly speaking, machine learning approaches do not benefit from prior information on the solution as described above, but ``learn'' the segmentation from large training datasets. The aim of a machine learning approach is to define a segmentation model $f_\theta:\mathcal{I} \longrightarrow \Sigma$ such that the segmentation of $I_0$ can be obtained as $I^*= f_\theta(I_0)$. The function $f_\theta$ is usually nonlinear and $\theta$ is a large vector of parameters. The learning phase selects $\theta $ in order to minimize a loss functional $\mathcal{L}$ that measures the accuracy of the predicted segmentation $f_\theta(I_0)$.

In \emph{supervised} machine learning, training data are available from databases of annotated segmentations, which provide a large number of pairs $(I_0,I^*) \in X \times Y \subset \mathcal{I} \times \Sigma$ \ ($X \times Y$ is named training set). The vector of parameters $\theta$ is obtained by minimizing a loss function plus a penalty term. For the sake of simplicity, we first consider a mean-square-error loss: 
\begin{equation}\label{eq:Cost}
   \theta^* = \arg \, \underset{\theta}{\min} \,\mathcal{L}(X,Y,\theta)= \arg \, \underset{\theta}{\min} \, \left( \sum_X \parallel f_\theta(I_0)- I^* \parallel^2 +
   \mathcal{P}_\theta (f_\theta(I_0)) \right).
\end{equation}
\noindent
Another widely used loss functions is the \textit{Binary Cross Entropy} (BCE) loss, which measures the difference in information content between the actual 
and the predicted image segmentation:
\[\mathcal{L}_{BCE}(f_\theta(I_0),I^*) = -f_\theta(I_0)\log(I^*) - (1-f_\theta(I_0))log(1-I^*) .\]
It is based on the Bernoulli distribution and works well with equal data distributions among classes. Some variants of BCE, such as the \textit{Weighted} BCE and the \textit{Balanced} CE are also used for tuning false negatives and false positives, respectively.
The \textit{Shape-aware} (Sa) loss 
calculates the average point-to-the-curve
Euclidean distance among points around the curve of the predicted segmentation, $u^*$, to the ground truth, $\bar{u}$, and use it as a coefficient to the cross-entropy ($CE$) loss function:
\[\mathcal{L}_{Sa}(f_\theta(I_0),I^*) = - \sum_{i \in \Sigma}
CE(\theta(I_0),I^*)- \sum_{i \in \Sigma}
i\, E_i \, CE(f_\theta(I_0),I^*),\]
where $\Sigma$ contains the
set of points  where the prediction curve does not match the ground-truth curve, and $E_i=d(u_i^*, \bar{u}_i)$.
The \textit{Dice} loss, based on the well-known Dice coefficient metric, is also widely used to measure the similarity between two segmentations, and is defined as
\[\mathcal{L}_{Dice}(f_\theta(I_0), I^*) = 1 - \frac{2f_\theta(I_0)I^*+1}{f_\theta(I_0)+I^*+1}.\]

In \emph{unsupervised} machine learning, the training set is not equipped with annotated segmentations and the goal is to train $f_\theta$ to recognize specific patterns or image features in the data. This approach is sometimes referred to as self-supervised learning \cite{bib:Doersch2015}, because the information is extracted from the data themselves rather than from a set of ``predictions'' (i.e., given segmentations). Then the fidelity term in~\eqref{eq:Cost} takes the form
\[ \sum_\mathcal{I} \parallel f_\theta(I_0)- \Phi(f_\theta(I_0))\parallel^2,\]
where $\Phi$ is the labeling operator defined in~\eqref{eq:labeling}.

In order to progressively extract higher-level features from the data, machine learning models use a multi-layer structure called \emph{neural network}, consisting of successive function compositions. The number of layers is the depth of the model, hence the terminology \emph{deep learning}. A neural network with $L$ layers is a function
\[ f_\theta : \mathcal{I} \times (H_1 \times \ldots \times H_L) \longrightarrow \Sigma, \quad f_\theta(I)=(f_L \circ f_{L-1} \circ ... \circ f_1)(I), \]
\noindent
where $f_i : \mathbb{R}^{d_{i-1}} \times H_i\longrightarrow \mathbb{R}^{d_{i}}$ are the activation functions (each depending on a component $\theta_i$ of $\theta$), \ $d_0 = d$ and $d_L = n$, with $n$ equal to the number of features. The adjective ``neural'' comes from the fact that those networks are loosely inspired by neuroscience.

\begin{figure}[ht!]
\medskip
\newcolumntype{C}{>{\centering\arraybackslash} m{.3\textwidth} }
\hspace*{-6pt}\begin{tabular}{CCC}
\includegraphics[width=.29\textwidth,height=.2\textwidth]{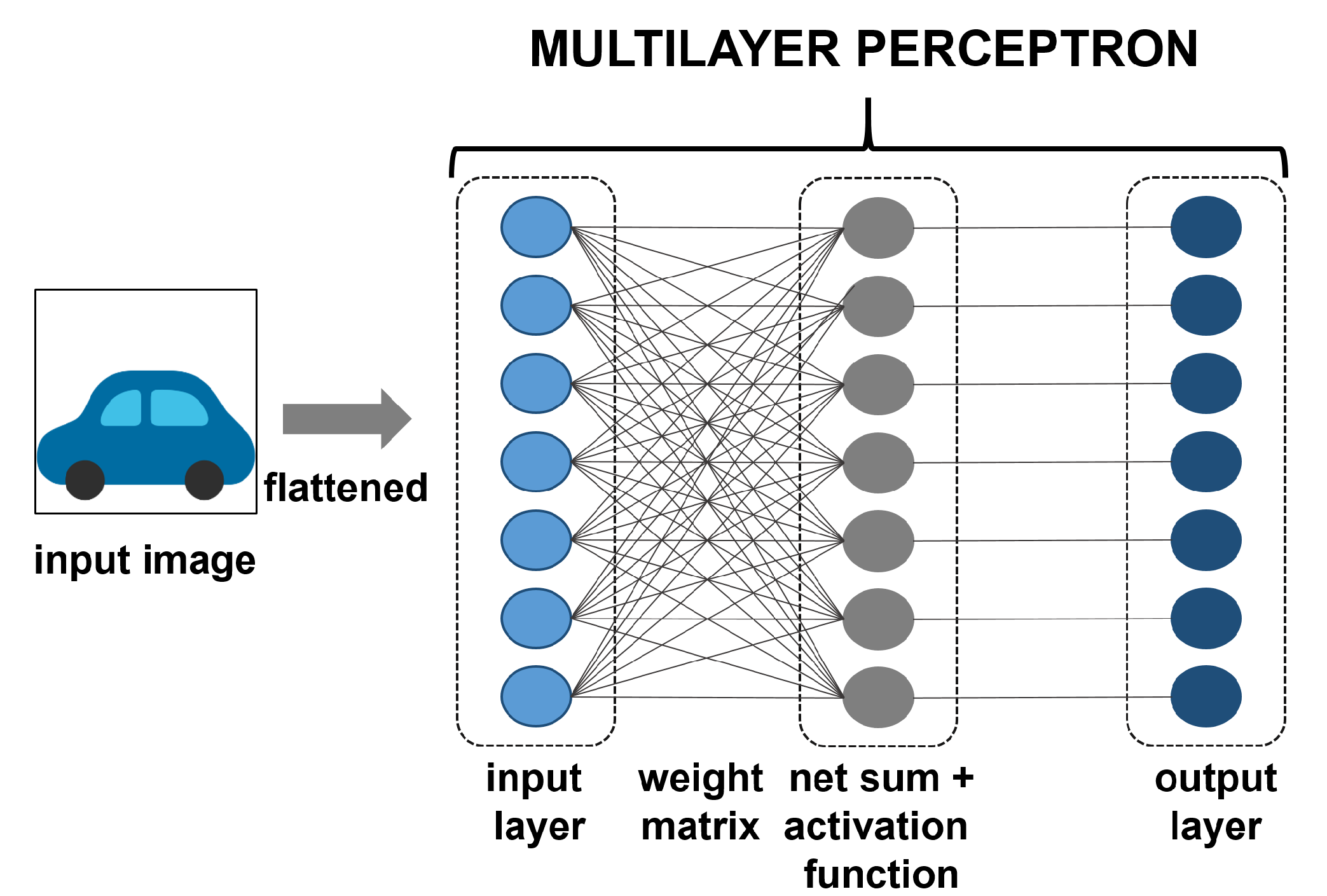} &
\includegraphics[width=.29\textwidth,height=.2\textwidth]{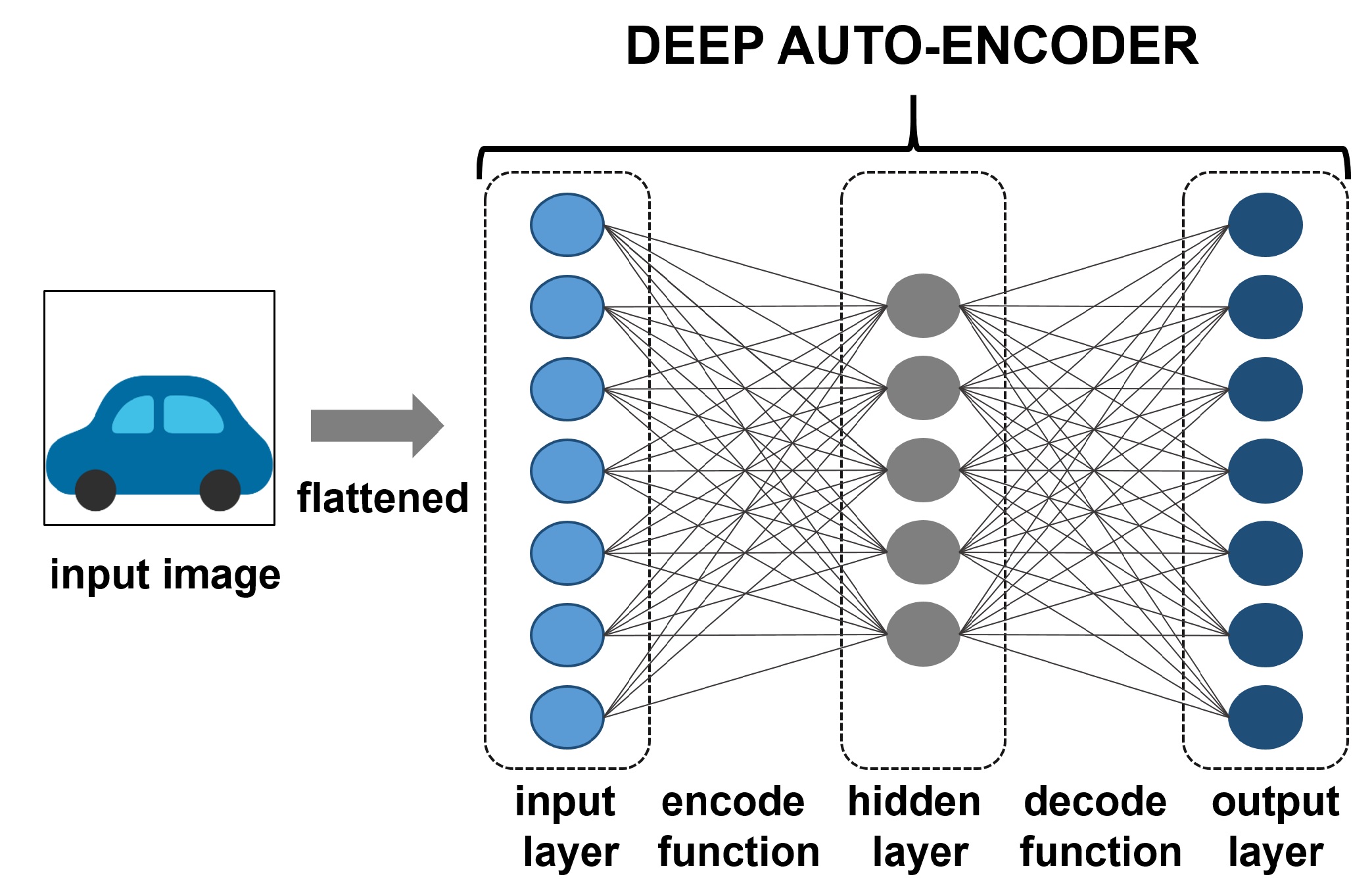} &
\includegraphics[width=.31\textwidth,height=.2\textwidth]{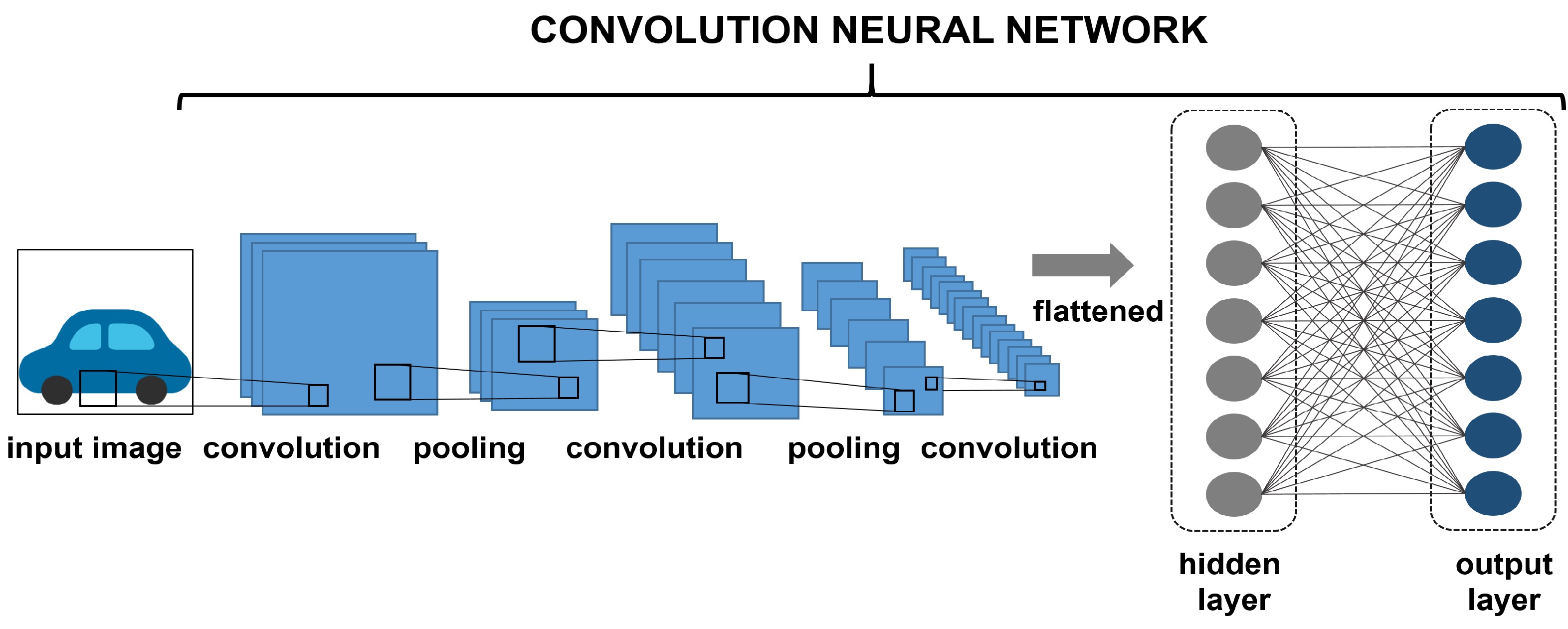}\\
\end{tabular}
\caption{Some neural network architectures used in image segmentation.\label{fig:DNN}}
\end{figure}

Neural network structures successfully used in image segmentation are the Multilayer Perceptron (MLP), the Deep Auto-Encoder (DAE) and the Convolutional Neural Network (CNN)~\cite{bib:Furat2019DL, bib:Minaee2020ImageSU, bib:Rizwan2020}.
Their basic schemes are shown in Figure~\ref{fig:DNN}.
%
The MLP is a neural network connecting multiple layers in a directed graph, which means that the signal path through the nodes only goes one way. Each node, apart from the input nodes, has a nonlinear activation function. An MLP uses backpropagation as a supervised learning technique.
The DAE network structure typically consists of $2L$ layer functions, where the first $L$ layers act as an encoding function with the input to each layer being of lower dimension than the input to the previous layer, and the remaining $L$ layers increase the size of their inputs until the final layer has the same dimension as the image input. The first $L$ layers are an MLP.
CNNs  divide the image into small areas and scan it one area at a time, to identify and extract features that are used to classify the image.
A CNN mainly consists of three layers: 
\begin{itemize}
    \item \emph{convolutional layer}: the image is analyzed a few pixels at a time to extract low-level features (edges, color, gradient orientation, etc.);
    \item \emph{nonlinear layer}: an element-wise activation function creates a feature map with probabilities that each feature belongs to the required class;
    \item \emph{pooling} or \emph{downsampling layer}: the amount of features and computations in the network is reduced, hence controlling overfitting.
\end{itemize}

\noindent Among well-known deep neural network architectures successfully used in image segmentation, we mention SegNet \cite{bib:SegNet2015}, U-Net \cite{bib:UNet2015},
and FCN \cite{bib:Long2015FCN}.


\section{Numerical techniques for segmentation models\label{sec:methods}}

The minimization in~\eqref{eq:generalE} is usually nontrivial and requires appropriate methods, taking into account the specific application.
In this section we provide a brief summary of numerical methods that can be applied to segmentation models. We consider two approaches: \emph{first discretize then optimize} and \emph{first optimize then discretize}. In the former, all the quantities in~\eqref{eq:generalE} are discretized a priori and then optimization methods are applied to the resulting minimization problem in $\mathbb{R}^n$. In the latter, we first write 
optimality conditions for~\eqref{eq:generalE}, which are generally partial differential equations (PDEs), and then solve those equations by suitable numerical methods, which discretize the equations. Finally, we also sketch some filtering techniques used in image segmentation, although they are not directly applied to the minimization problem~\eqref{eq:generalE}. This is motivated by their use in some segmentation approaches, such as those based on deep learning.

For the sake of simplicity, here we consider $S = I$ (i.e., we neglect $u$ in the segmentation $S = (I, u)$).
For 2D images ($d=2$)
we denote by $\Omega_{n_x,n_y}$ the discretization of $\Omega$ consisting of a grid of $ n_x \times n_y$ pixels,
\[ \Omega_{n_x,n_y} = \{ (i,j) : \, i=1,...,n_x, \; j=1,...,n_y\}. \]
\noindent
We also identify each pixel with its center and denote by $S_{i,j}$ the value of $S$ in $(i,j)$. Finally, we consider the forward and backward difference operators defined as follows:
\[
\begin{array}{ll}
\displaystyle
D^{+x} I_{i,j} = I_{i+1,j}-I_{i,j}  , &
\displaystyle
D^{-x} I_{i,j} = I_{i,j}-I_{i-1,j}  , \\[2mm]
\displaystyle
D^{+y} I_{i,j} = I_{i,j+1}-I_{i,j} , &
\displaystyle
D^{-y} I_{i,j} = I_{i,j}-I_{i,j-1} ,  \\[2mm]
\end{array}
\]
\noindent
where we assume
\[
\begin{array}{ll}
\displaystyle I_{i-1,j} = I_{i,j} \;\; \mbox{for } i=1,  & \displaystyle I_{i,j-1} = I_{i,j} \;\; \mbox{for } j=1,  \\[1mm]
\displaystyle I_{i+1,j} = I_{i,j} \;\; \mbox{for } i=n_x, & \displaystyle  I_{i,j+1} = I_{i,j} \;\; \mbox{for } j=n_y,
\end{array}
\]
\noindent
i.e., we define by replication the values of $I$ with indices outside~$\Omega_{n_x,n_y}$. Likewise, for 3D images the discretization of the image domain consists of a grid of $n_x \times n_y \times n_z$ voxels,
\[ \Omega_{n_x,n_y,n_z} = \{ (i,j,k) : \, i=1,...,n_x, \; j=1,...,n_y, \; k=1,...,n_z\}, \]
\noindent
and the forward and backward difference operators are defined as follows:
\[
\begin{array}{ll}
\displaystyle
D^{+x} I_{i,j,k} = I_{i+1,j,k}-I_{i,j,k}  , &
\displaystyle
D^{-x} I_{i,j,k} = I_{i,j,k}-I_{i-1,j,k},\\[2mm]
\displaystyle
D^{+y} I_{i,j,k} = I_{i,j+1,k}-I_{i,j,k} , &
\displaystyle
D^{-y} I_{i,j,k} = I_{i,j,k}-I_{i,j-1,k},\\[2mm]
\displaystyle
D^{+z} I_{i,j,k} = I_{i,j,k+1}-I_{i,j,k}, &
\displaystyle
D^{-z} I_{i,j,k} = I_{i,j,k}-I_{i,j,k-1}.\\[2mm]
\end{array}
\]
\noindent
For simplicity, henceforth we consider $d=2$.


\subsection{First discretize then optimize}

Numerical optimization offers a large variety of methods to compute the segmentation by solving the minimization problem coming from a discretization of~\eqref{eq:generalE}, possibly subject to constraints that can drive the segmentation towards particular features. The choice of the optimization method depends on the properties of the objective function and/or the constraints.

Roughly speaking, at iteration $k$, optimization methods for nonlinear problems generate a function $\widetilde{E}(I;I_k)$ that approximates the discretized objective function $E$ around $I_k$, and minimize it to obtain the next iterate (see, e.g., \cite{bib:FontoulakisGondzio2016}). For example, given $I_k$, the $(k+1)$-st iteration may be written as
\[  \begin{array}{l}
   \mbox{Define } \widetilde{E}(I;I_k)  \mbox{ that approximates } E(I; I_0) \\[2pt]
   \widetilde{I}_{k+1} = \arg \underset{I}{\min} \; \widetilde{E}(I;I_k) \\
   I_{k+1} = I_k + \alpha_k (\widetilde{I}_{k+1}-I_k)
\end{array} \]
\noindent
where the step length $\alpha_k$ satisfies some criterion.

``Classical'' optimization techniques, such as gradient or Newton-type methods, 
require regularity assumptions on the objective function (and the constraints, if any). However, many segmentation models are modeled as non-smooth optimization problems. There are two main approaches to deal with non-differentiability: smoothing and non-smoothing~\cite{bib:Antonelli18}. The former formulates the problem as a suitable smooth one and applies the aforementioned classical optimization methods. The latter does not modify the mathematical model, and thus uses methods not requiring smoothness. For the purpose of illustration, here we focus on~\eqref{eq:CENmodel}, where non-smoothness comes from a discretization of the TV term.

A regularized discrete TV may be obtained as follows:
\[ \int_\Omega \vert \nabla I \vert \, dx \approx \sum_{i,j} \sqrt{(D^{+x} I_{i,j})^2+(D^{+y} I_{i,j})^2+\epsilon}, \]
\noindent
where $\epsilon > 0$ is ``suitably small'', but other regularized versions may be considered, e.g., based on Huber-like functions~\cite{bib:WeissBlancFeraudAubert2009}. In this case, gradient and higher-order methods~\cite{bib:SPG2000,bib:BonettiniZanellaZanni2009,bib:Antonelli2016SPGA,bib:diserafino18AMC,bib:diserafino20AMC} 
can be used efficiently. Another way of introducing smoothness consists in splitting the variables into their positive and negative parts (thus doubling the number of unknowns) and introducing new constraints, and then applying first- or higher-order methods for smooth problems, such as in~\cite{bib:FigueiredoNowakWright2007,bib:FontoulakisGondzioZhlobich2014,bib:desimone2022}.

Non-smoothing approaches avoid regularization of the non-smooth terms in the optimization problem. This is the case, for example, of methods based on forward-backward splitting techniques, such as proximal-gradient methods~\cite{bib:ParikhBoyd2014,bib:Bonettini2016}, and the forward-backward Expectation Maximization (EM) method in~\cite{bib:Sawatzsky2008}. ADMM and split Bregman methods do not use smooth approximations too~\cite{bib:Boyd2011,bib:FigueiredoBioucas-Dias2010,bib:GoldsteinBressonOsher2010,bib:Setzer:2011,bib:desimone2020,bib:Antonelli2020Adaptive}. The success of these approaches is based also on the availability in closed (and cheap) form of the proximal operator of the $\ell_1$ norm by means of the well-known soft-thresholding, defined as
$$
[{\mathcal S}(x,\gamma)]_{i,j}= \mathrm{sign}(x_{i,j})\cdot\max\big(\vert x_{i,j}\vert-\gamma, 0\big),
$$
with $\gamma > 0$. The difficulties associated with the non-differentiability of the TV functional may be also overcome by reformulating the minimization problem as a saddle-point problem and solving it by a primal-dual algorithm such as the Chambolle-Pock one~\cite{bib:ChambollePock2011,bib:MalitskyPock2018}.

EM algorithms~\cite{bib:EMalgo} are also widely used to solve statistical models. They are based on the idea of splitting the (negative) log-likelihood into two terms and alternating between the computation of the expectation and its minimization.


Finally, stochastic versions of the previous methods are used in segmentation with deep learning, to limit the computational cost. The idea is to use only random samples of the data at each iteration, to estimate first-order and possibly second-order information according to the loss function, with the aim of significantly reducing the computation and hence the time~\cite{bib:Guanci2018,bib:stochastic19}.


\subsection{First optimize then discretize}

Reducing imaging problems to PDEs is many years old, because of the availability of a large amount of methods and software for solving PDEs. PDE-based methods have been introduced in different ways, such as the Perona-Malik filtering~\cite{bib:PeronaMalik}, directly based on properties of the PDE~\cite{bib:witkin}, and the axiomatic scale space theory~\cite{bib:scalespace,bib:AlvarezMorel2}.

In a variational approach, one derives the first-order optimality conditions via smoothing regularization, if it is needed. Let us consider, for example,
the level-set formulation of the Chan-Vese model~\eqref{eq:CVmodel}, where $I$ is represented by a function $\phi$ such that $\phi(x) = 0$ provides the curve separating two regions of $I$ (when $I=I^*$ the two regions identify the segmentation). Keeping $c_{in}$ and $c_{out}$ fixed and writing the Euler-Lagrange equations in a gradient-flow approach, we get
\begin{equation}\label{eq:cv_gradientflow}
\begin{array}{llll}
  \displaystyle \frac{\partial\phi}{\partial t} (t,x) & = & \displaystyle \delta_\varepsilon (\phi) \left( \lambda \, \mbox{div} \left( \frac{\nabla \phi}{\vert \nabla \phi \vert} \right) - (c_{in} - I_0)^2 + (c_{out} - I_0)^2 \right) & \mbox{in } (0, +\infty) \times \Omega , \\[7pt]
  \phi(0,x) & = & \phi_0(x) & \mbox{in } \Omega , \\[3pt]
  \displaystyle \frac{\delta_\varepsilon (\phi)}{\vert \nabla \phi \vert}\frac{\partial \phi}{\partial N} & = & 0 & \mbox{on } \partial\Omega ,
\end{array}
\end{equation}
\noindent
where $\delta_\varepsilon$ is a regularized version of the Dirac measure, $\phi_0$ is the initial-level function, and $N$ is the exterior normal to the boundary $\partial\Omega$~\cite{bib:ChanVese2001}.

Finite-difference schemes 
are popular methods for the numerical solution of~\eqref{eq:cv_gradientflow}.
Of course, the discretization used in image segmentation must take into account the nature and the properties of the operators involved in the model. For example, edge preserving is similar to shock capturing in computational fluid dynamics, and hence finite-difference schemes based on hyperbolic conservation laws can be used~\cite{bib:Sethian99}.
Just to give an example, the level-set equation
\[ \frac{\partial \phi}{\partial t} = F \vert \nabla \phi \vert \]
in Subsection~\ref{subsec:edge-based}
can be solved by using an upwind numerical scheme:
\[ \phi^{n+1} = \Psi(\phi^n), \quad \Psi(\phi^n_{i,j}) = \phi^n_{i,j} - \Delta t (\max(F,0) \nabla^+ \phi^n_{i,j} + \min(F,0)\nabla^- \phi^n_{i,j}) , \]
\noindent
where
\[
\begin{array}{lll}
\displaystyle \nabla^+ \phi^n_{i,j} & = & \displaystyle \left( \max\left( \max(D^{-x} \phi^n_{i,j}, \, 0)^2, -\min(D^{+x} \phi^n_{i,j}, \, 0)^2 \right) \right. \\
                                                    & + & \displaystyle \left. \max\left( \max(D^{-y} \phi^n_{i,j}, \, 0)^2, -\min(D^{+y} \phi^n_{i,j}, \, 0)^2 \right) \right)^{1/2} ,\\[7pt]
\displaystyle \nabla^- \phi^n_{i,j} & = & \displaystyle \left( \max\left( \max(D^{+x} \phi^n_{i,j}, \, 0)^2, -\min(D^{-x} \phi^n_{i,j}, \, 0)^2 \right) \right.\\
                                                    & + & \displaystyle \left. \max\left ( \max(D^{+y} \phi^n_{i,j}, \, 0)^2, -\min(D^{-y} \phi^n_{i,j}, \, 0)^2 \right) \right)^{1/2} .
\end{array}
\]


\subsection{Filters}

Discrete filters are often used in image segmentation, e.g., in machine learning approaches. A digital filter can be represented as an operator
\[ L: I \in \mathcal{I} \longrightarrow \widetilde{I} \in \mathcal{I}, \quad \widetilde{I}_{ij} = L [I; W_{ij}], \]
where $W_{ij} \subset \Omega_{n_x,n_y} $. A popular discrete filter in image segmentation is the \emph{convolution} filter, defined by
\begin{equation}\label{eq:convolution}
\widetilde{I}_{i,j} = L_{a,b}[I; W_{i,j}] = \sum_{s=-a}^a \sum_{t=-b}^b h_{s,t} I_{i-s,j-t} ,
\end{equation}
\noindent
with $a$ and $b$ positive integers such that $a \leq \frac{n_x-1}{2}$ and $b \leq \frac{n_y-1}{2}$, $W_{i,j} = \{ (s,t) : s = -a, \ldots, a, \; t = -b, \ldots, b \} $, and $h_{s,t} \in \mathbb{R}$. The matrix $H = (H_{i,j}) = (h_{-a+i,-b+j}) \in \mathbb{R}^{(2a+1)\times (2b+1)}$ is called kernel matrix and depends on the features we want to extract from the image. Common choices of $a$ and $b$ are $a=b=3$ and $a=b=5$.

Edge-detection kernels are frequently used in image segmentation, especially in CNNs. For example, the first layer of a CNN is often responsible for capturing low-level features such as edges, color, and gradient orientation. In general, the choice of $H$ determines the type of features to be extracted. The kernel matrix
\[ H= \left ( \begin{array}{lll}
     1 & \;\, 0 &-1  \\
     1 & \;\, 0 &-1\\
     1 & \;\, 0 &-1
     \end{array}
\right ) \]
\noindent
is a vertical edge-detection kernel~\cite{bib:edgedetector2019}. Another example is the Sobel operator, used to create an image emphasizing the edges~\cite{bib:kanopoulos1988design}. It allows us to obtain either the gradient amplitude or the gradient direction of the image intensity at each point, by convolving the image with the kernel matrices
\[ H^x_S=
   \left ( \begin{array}{rrr}
    1 & \;\, 0 & -1 \\
    2 & \;\, 0 & -2 \\
    1 & \;\, 0 & -1
   \end{array}
   \right ), \quad
     H^y_S= \left ( \begin{array}{rrr}
    1 & 2 & 1 \\
    0 & 0 & 0 \\
   -1 &-2 &-1
   \end{array} \right ).
\]
\noindent
The gradient magnitude, $G$, and the angle of orientation of the edges, $\theta$, are given by
\[
\lvert G_{i,j} \rvert = \sqrt{(H^x_S*I)_{i,j}^2+(H^y_S*I)_{i,j}^2}, \quad \theta_{i,j}=\arctan((H^y_S*I)_{i,j}/(H^x_S*I)_{i,j}).
\]
\noindent
A padding process is commonly used to preserve the dimension of the image after the convolution. It usually consists in the replication or reflection of the pixel values at the image border, or in adding an average gray or even zeros symmetrically around the border of the image. A pooling layer is usually inserted between two successive convolution layers, which is obtained by applying basic functions, such as max and mean, in a small window.


\section{Conclusion}\label{sec:conclusions}

We presented a view of image segmentation, focusing on simple computational models and attempting to put them into a coherent framework where the inclusion of a-priori information about the solution is obtained by using penalty terms. We first introduced image segmentation and then outlined basic edge-based, region-based, statistical and machine learning models. We also sketched some numerical methods that can be employed to compute solutions to the models. We believe that our view of models and methods for image segmentation, although very far from being exhaustive, can help the readers understand much modern and sophisticated segmentation techniques, as well as select computational tools for their problems.


\backmatter





\bmhead{Acknowledgments}
This work was partially supported by the Istituto Nazionale di Alta Matematica - Gruppo Nazionale per il Calcolo Scientifico (INdAM-GNCS), by the Italian Ministry of University and Research under grant no. PON03PE\_00060\_5, and by the VALERE Program of the University of Campania ``L. Vanvitelli''. We would like to thank Giuseppe Trerotola (ICAR-CNR) for his technical support.

\bibliography{biblio_ReviewSeg}


\end{document}